\documentclass[12pt]{article}
\usepackage{amsmath,amssymb,hhline}
\usepackage[english,russian]{babel}

\begin{document}
\noindent{УДК 517.5} \vskip 0.1cm

\noindent{\large\bf R.~M.~Trigub}
\vskip 0.3cm

\noindent{\bf ON APPROXIMATION OF FUNCTIONS BY POLYNOMIALS \\
AND BY ENTIRE FUNCTIONS OF EXPONENTIAL TYPE} \vskip 0.2cm

\noindent{\large\bf Р.~М.~Тригуб}
\vskip 0.3cm

\noindent{\bf О ПРИБЛИЖЕНИИ ФУНКЦИЙ ПОЛИНОМАМИ И ЦЕЛЫМИ \\ФУНКЦИЯМИ
ЭКСПОНЕНЦИАЛЬНОГО ТИПА} \vskip 0.2cm

\noindent{\footnotesize A brief overview of publications in approximation
theory of functions known to the author and connected with scientific
publications by V.~K.~Dzyadyk (1919--1998).}

\vskip 1mm


 \linespread{1.2}{
Владислав Кириллович Дзядык  закончил Днепропетровский университет (ДГУ) в 1951~г. и в том же году вышла в Докл. АН СССР его заметка. Уже
в 1953 г. его статья, в которой решена проблема Фавара для $r\in(0,1)$ появилась в Изв.~АН СССР (см. ниже в п.~1). 
В настоящей статье излагается краткий обзор теорем по теории приближений, примыкающих к научным результатам В.К. См. также в п.~4 некоторые нерешенные вопросы.

Кстати, в 2019 исполняется 100 лет теории приближений, как ветви
математического анализа, если считать началом выход книги
\cite{ValleePoussin} Валле--Пуссена, в которой собраны теоремы
П.~Л.~Чебышева, К.~Вейерштрасса (K. Weierstrass), А.~А.~Маркова, А.~Лебега (H. Lebesgue), Л.~Фейера (L. Fej\'er),
  Д.~Джексона (D. Jackson),  С.~Н.~Бернштейна и самого автора.

\bigskip

\textbf{1. Наилучшее приближение класса $W^r(\mathbb{T})$
периодических функций полиномами и класса $W^r(\mathbb{R})$ целыми
функциями экспоненциального типа (ЦФЭТ).}

 Это классы функций с ограниченным единицей модулем производной $f^{(r)}$ (при $r$ нецелом --- производная по Вейлю, а при $r\in\mathbb{N}$ --- это функции, у
которых $f^{(r-1)}\in {\rm AC_{loc}}$, а $\big|f^{(r)}(x)\big|\leq1$ почти всюду).

Если функция $2\pi$--периодическая и принадлежит $L_p(\mathbb{T})$,
$\mathbb{T}=[-\pi,\pi]$, $p\in[1,+\infty]$ или $C(\mathbb{T})$, то
ее ряд Фурье будем писать в виде $\big(e_k=e^{ikx},\
k\in\mathbb{Z}\big)$
\begin{equation*}
    f\sim\sum\limits_{k\in\mathbb{Z}}\widehat{f}_ke_k,\qquad
    \widehat{f}_k=\frac{1}{2\pi}\int\limits_{\mathbb{T}}f(t)e^{-ikt}dt.
\end{equation*}

Известно, что ($b_r$ --- ядро Бернулли)
\begin{equation*}
    f(x)=\frac{1}{\pi}\int\limits_{-\pi}^{\pi}f^{(r)}(x-t)b_r(t)dt,\qquad
    b_r(t)=\sum\limits_{k=1}^\infty\frac{\cos\big(kt-\frac{r\pi}{2}\big)}{k^r}.
\end{equation*}

При натуральном $r$ $b_r$ является периодическим интегралом  от
$\displaystyle b_1(t)=\frac{1}{2}\Big(1-\frac{|t|}{\pi}\Big){\rm sign}~t$ ($|t|\leq\pi$).

Г.~Бор (H. Bohr) в 1935 г. отметил без доказательства, что если $r=1$ и
$\widehat{f_0}=0$, то $$\|f'\|_\infty\leq\frac{\pi}{2}\|f\|_\infty,$$
и неравенство точное. Затем С.~Н.~Бернштейн (1935)
доказал, что при $r\in\mathbb{N}$ и $\widehat{f_0}=0$\ $\|f^{(r)}\|_\infty\leq K_r\|f\|_\infty$, где
\begin{equation*}
    K_r=\frac{4}{\pi}\sum\limits_{\nu=0}^\infty\frac{(-1)^{\nu(r+1)}}{(2\nu+1)^{r+1}}.
\end{equation*}

Ж.~Фавар (J. Favard) \cite{Favard36}, цитируя статью Бернштейна, определил
наилучшее приближение в $L_1(\mathbb{T})$ ядра Бернулли
тригонометрическими полиномами $\tau_n$ порядка не выше $n$
$\big(\tau_n=\sum\limits_{|k|\leq n}c_ke_k\big):$
\begin{equation*}
    E_n^T(b_r)_1=\min\limits_{\tau_n}\int\limits_{\mathbb{T}}\big|b_r(t)-\tau_n(t)\big|dt=\frac{\pi
    K_r}{(n+1)^r}.
\end{equation*}

Затем одновременно и независимо Фавар \cite{Favard37} и Ахиезер-Крейн \cite{Akhiezer_Krein} вывели из этого соотношения, что
\begin{equation*}
    \sup\limits_{f\in W^r(\mathbb{T})}E_n^T(f)_\infty=\sup\limits_{f\in
    W^r(\mathbb{T})}\min\limits_{\tau_n}\sup\limits_x\big|f(x)-\tau_n(x)\big|=\frac{K_r}{(n+1)^r}.
\end{equation*}

Экстремальная функция для класса равна (сплайн Эйлера)
\begin{equation*}
\varphi_r(x)=\frac{4}{\pi}\sum\limits_{k=0}^\infty\frac{\sin\big((2k+1)x-\frac{r\pi}{2}\big)}{(2k+1)^{r+1}}.
\end{equation*}

Это $r$-й периодический интеграл от $\varphi_0(x)={\rm sign}\sin x$.
В \cite{Akhiezer_Krein} еще изучен тот же вопрос для класса
$\widetilde{W_r}(\mathbb{T})$
$\big(\big\|\widetilde{f}^{(r)}\big\|_\infty\leq 1$, $\widetilde{f}$
--- тригонометрически сопряженная к $f$ функция$\big)$ с другой константой $\widetilde{K_r}$.

Бернштейн назвал, по сути, $K_r$ константой Бернулли при нечетном
$r$ и константой Эйлера при чeтном $r$. Так что нельзя, я считаю,
называть $K_r$ константой Фавара. Статья Бернштейна обычно не
цитируется, так как в ней в доказательстве при четном $r$ была
ошибка (исправлена позже в \cite[с.~171]{Bernstein}).

Вообще, задача о наилучшем приближении класса функций (а не
индивидуальных функций) инициирована А.~Н.~Колмогоровым (1935,
\cite{Kolmogorov}). Ею успешно занимался С.~М.~Никольский (1905-2012) и многие другие математики.

Фавар (1937) поставил задачу о наилучшем приближении класса $W^r(\mathbb{T})$ при нецелом $r>0$. 
Решение ее существенно труднее, так как ядро Бернулли $b_r$ не является ни четным, ни нечетным.
Ответ при $r\in(0,1)$ получил В.К. Дзядык \cite{Dzjadyk53}. Вслед за ним этой задачей занимались С.~Б.~Стечкин и Сунь Юн--Шен (Sun Yung-sheng).
Но завершил ее решение для любого $r>0$ В.К. Дзядык, применяя новый метод с абсолютно монотонными функциями \cite{Dzjadyk74} (в метрике $L_\infty$ и $L_1$).

В 1957 В.К. защитил кандидатскую диссертацию в ДГУ и вернулся в Луцкий пединститут. В 1959 в Ленинграде состоялась  первая Всесоюзная конференция
по конструктивной теории функций, на которую из ДГУ поехали во главе с профессором А.Ф. Тиманом девять человек, из которых пять были еще аспирантами,
включая автора этих строк. На этой конференции, в частности, В.~М.~Тихомиров (Москва) сообщил, что тригонометрические полиномы $\tau_n$ осуществляют 
наилучшее приближение класса $W^r(\mathbb{T})$ среди всех подпространств размерности $2n+1$ $\big($поперечник Колмогорова
$\displaystyle
d_{2n+1}\big(W^r(\mathbb{T})\big)=\frac{K_r}{(n+1)^r}\big)$.\vskip
1mm

Случай метрики $L_p$, $p\in(1,+\infty)$, особый, и некоторые неравенства сразу следуют, как оказалось, из случая $p=\infty$. Но при этом может быть потеряна точность.

Введем, например, семейство норм $(N>0)$
\begin{equation*}
    \|f\|_{p,N}=\sup\limits_{x\in\mathbb{R}}\Big(\int\limits_{x-N}^{x+N}|f(t)|^pdt\Big)^{\frac{1}{p}}.
\end{equation*}

Если $f\in L_p(\mathbb{R})$, то
$\|f\|_p=\lim\limits_{N\rightarrow\infty}\|f\|_{p,N}$, а если $f$
--- $2\pi$--периодическая и $f\in L_p(\mathbb{T})$, то
$\|f\|_p=\|f\|_{p,\pi}$. Тем самым еще соединяем периодический и непериодический случаи.

\textbf{Теорема} (\cite[\textbf{1.2.7}]{Trigub_Belinsky})
\emph{Пусть $E$ -- линейное множество ограниченных и равномерно непрерывных на
$\mathbb{R}$ функций, замкнутое относительно равномерной сходимости,
а $$A:~E~\mapsto~E$$  линейный оператор и $\|Af\|_\infty\leq
a\|f\|_\infty$. Если дополнительно вместе с $f$ и $f^t=f(\cdot+t)\in
E$ при $t\in\mathbb{R}$ и оператор $A$ коммутирует со сдвигом, то
для любого $p\in[1,+\infty)$, $N>0$ и $f\in E$}
\begin{equation*}
    \|Af\|_{p,N}\leq a\|f\|_{p,N}.
\end{equation*}

\medskip

\textbf{Следствия.}

\emph{I. Точное неравенство Бернштейна для полиномов и ЦФЭТ по
$L_p$-норме следуют из случая $p=\infty$ для полиномов. А из
неравенства Бернштейна для ЦФЭТ в $L_2$,кстати, следует классическая
теорема Винера--Пэли. См. \cite[с.~89]{Trigub_Belinsky}.}

\emph{II. Теорема типа Джексона (см. в п.~2) с линейными операторами
в метрике $L_p$ следует из случая $p=\infty$.}

Непериодические функции на прямой приближают, следуя Бернштейну,
ЦФЭТ. Любая ЦФЭТ, ограниченная на $\mathbb{R}$, является пределом
последовательности тригонометрических полиномов, сходящейся
равномерно на любом отрезке $\mathbb{R}$ (Б.~М.~Левитан). См.,
например, \cite{Akhiezer} -- \cite{DeVore}, \cite[\textbf{4.2.8}]{Trigub_Belinsky}.
И в следующей лемме содержится переход от периодического случая к непериодическому.

Обозначим через $g_\sigma$, $\sigma>0$, --- ограниченную на
$\mathbb{R}$ ЦФЭТ $\leq\sigma$. Если $g_\sigma$ имеет период $2\pi$,
то $g_\sigma=\tau_{[\sigma]}$ (полином).

Положим
\begin{equation*}
    A_\sigma(f)_\infty=\inf\limits_{g\sigma}\sup\limits_{\mathbb{R}}\big|f(x)-g_\sigma(x)\big|.
\end{equation*}

\textbf{Лемма} (\cite[\textbf{5.5.9}]{Trigub_Belinsky}) \emph{Пусть
$W$ -- множество ограниченных непрерывных функций на $\mathbb{R}$,
замкнутое относительно сходимости почти всюду и преобразования
подобия (если $f\in W$, то и $f_\lambda\in W$ при $\lambda>0$, где
$f_\lambda(x)=f(\lambda x)$). Пусть еще любая функция $f\in W$
является пределом последовательности периодических функций из $W$.
Если, кроме того, имеется функция $\varepsilon:
W\times(0,+\infty)\mapsto(0,+\infty)$ с двумя свойствами:
$\varepsilon(f_\lambda,h)=\varepsilon(f,\lambda h)$ $(h>0)$ и из
того, что $f_m\rightarrow f$ при $m\rightarrow\infty$, следует, что
$\varepsilon(f_m,h)\rightarrow\varepsilon(f,h)$, то из неравенства
$\displaystyle
E_n^T(f)_\infty\leq\varepsilon\Big(f;\frac{1}{n}\Big)$
$(n\in\mathbb{N})$ для $2\pi$--периодических функций из $W$ следует,
что для любой функции $f\in W$ при $\sigma>0$}
\begin{equation*}
    A_\sigma(f)\leq\varepsilon\Big(f;\frac{1}{\sigma}\Big).
\end{equation*}

\medskip

\textbf{Следствия из периодического случая.}\vskip 1mm

I. $\displaystyle
A_\sigma\big(W^r(\mathbb{R})\big)_\infty=\frac{K_r}{\sigma^r}\quad
(\sigma>0,\ r\in\mathbb{N}).$\vskip 1mm

II. $\displaystyle
A_\sigma(C_\omega)_\infty=\frac{1}{2}\omega\Big(\frac{\pi}{\sigma}\Big)\quad
(\sigma>0).$\vskip 1mm

Прямое доказательство первого соотношения см. в \cite{Akhiezer}. А
используя теорему В.К.Дзяды\-ка  получаем, что соотношение I верно при любом $r>0$ (с другой константой вместо $K_r$ при нецелом $r$).

Второе равенство, а здесь речь идет о наилучшем приближении класса
непрерывных функций с заданной выпуклой мажорантой $\omega$ модулей
непрерывности, для периодических функций доказано Н.~П.~Корнейчуком (1920-2003)
(старая задача, доложенная Н.~П. впервые в марте 1961 на семинаре
А.~Ф.~Тимана в ДГУ). В приведенном виде -- это теорема В.К. Дзядыка (см. \cite{Dzjadyk75}).

О науке в Днепропетровске см. в воспоминаниях С.~М.~Никольского (МИРАН, М., 2003, 160~с.; доступны на сайте Mathnet.ru).

Отметим ещё, что уже после работ А.~А.~Маркова и С.~М.~Никольского
критерий наилучшего приближения в $L_1$ принял следующий окончательный вид:

\textbf{Лемма} (\cite[\textbf{5.2.5}]{Trigub_Belinsky}) \emph{Пусть
$E$ -- подпространство $L_1$ и $f\in L_1\setminus E$. Для того чтобы}
\begin{equation*}
    {\rm dist}~(f,E)=\min\limits_{g\in E}\|f-g\|_{L_1}=\big\|f-g^*\big\|_{L_1}
\end{equation*}
\emph{необходимо и достаточно, чтобы существовала функция $h$ такая,
что $$|h|\leq1, h(f-g^*)=\big|f-g^*\big|$$ почти всюду и $\int
hg=0$ для всех $g\in E$. Кроме того, ${\rm dist}~(f,E)=\int hg$.
Если дополнительно $f-g^*\neq0$ почти всюду, то $h=\overline{{\rm
sign}~(f-g^*)}$.}
\\

\bigskip

\textbf{2. Приближение функций многочленами.}

\medskip

 Переходим к прямым и обратным теоремам теории приближений функций тригонометрическими и алгебраическими полиномами.

Бернштейн доказал, что для $2\pi$--периодических функций из ${\rm
Lip}~\alpha$ при $\alpha\in(0,1)$ и средних арифметических сумм Фурье (суммы Фейера $\sigma_n$)
\begin{equation*}
    \max\limits_x\big|f(x)-\sigma_n(x)\big|=O\Big(\frac{1}{n^\alpha}\Big),
\end{equation*}
а при $\alpha=1$ --- $\displaystyle O\Big(\frac{\ln n}{n}\Big)$ (см. \cite[с.~89]{Bernstein} или \cite[с.~236]{DeVore}).

Д.~Джексон независимо, заменив в $\sigma_n$ квадрат ядра Дирихле $D_n^2$  на $D_n^4$ с соответствующей нормировкой получил для 
${\rm Lip}~1$, а это главный случай, порядок приближения $\displaystyle O\Big(\frac{1}{n}\Big)$ (см.  \cite{Akhiezer} -- \cite{DeVore}). 
Отсюда уже легко вывести, используя функцию типа Стеклова, общую теорему Джексона--Стечкина с модулем гладкости
$\omega_r$ любого порядка (см., например, \cite{Trigub_Belinsky}).

Есть другой метод доказательства -- метод мультипликаторов Фурье
(см. \cite[ гл.~8]{Trigub_Belinsky}). Как показано в \cite{Trigub2009_East} этим только методом, точный порядок приближения $\displaystyle
\omega_2\Big(f;\frac{1}{n}\Big)$ один и тот же при $D_n^3$ и $D_n^4$.

Автором найден точный порядок приближения для классических
методов суммирования рядов Фурье (см. \cite{Dzjadyk77},
\cite{Trigub2013}). В случае функций многих переменных пришлось
вводить специальные модули гладкости (\cite{Trigub80}).

Если функция задана на отрезке ($[-1,1]$, например), то после
стандартной замены $x=\cos t$ функция становится $2\pi$--периодической
и четной (той же гладкости, в общем). А многочлен (алгебраический полином) $p_n$ степени не выше $n$ 
после замены $x=\cos t$ становится полиномом $\tau_n$. Получаем сразу, например , что
\begin{equation*}
    E_n(f)_\infty=\min\limits_{p_n}\max\limits_{[-1,1]}\big|f(x)-p_n(x)\big|=O\Big(\frac{1}{n^\alpha}\Big)\qquad
    (\alpha\in(0,1))
\end{equation*}
тогда и только тогда, когда $f(\cos t)\in {\rm Lip}~\alpha$ или для любых $x_1$ и $x_2$ из $[-1,1]$
\begin{equation*}
    \big|f(x_1)-f(x_2)\big|\leq\gamma\Big(\frac{|x_1-x_2|}{|x_1-x_2|+\sqrt{1-x_1^2}+\sqrt{1-x_2^2}}\Big)^\alpha.
\end{equation*}

Существенно более общие подобные теоремы о приближении многочленами
аналитических функций на компактах с углами на границе получены в статье \cite{Dzjadyk_Alibekov}.

С.~М.~Никольский (1946) доказал, что для любой функции $f\in W^1[-1,1]$ найдется
последовательность многочленов $p_n$ такая, что при $x\in[-1,1]$
\begin{equation*}
    \big|f(x)-p_n(x)\big|\leq\frac{\pi}{2}\cdot\frac{\sqrt{1-x^2}}{n}+O\Big(\frac{\ln
    n}{n^2}\Big).
\end{equation*}

Здесь важны поточечная оценка (у концов отрезка приближение почти в
два раза лучше) и наименьшая возможная константа $\displaystyle \frac{\pi}{2}=K_1$.

А.~Ф.~Тиман (1951) получил общую теорему с модулем непрерывности $\omega\big(f^{(r)};h\big)$

\noindent
$\displaystyle\big(\delta_n(x)=\frac{\sqrt{1-x^2}}{n}+\frac{1}{n^2}\big)$:
\begin{equation*}
    \big|f(x)-p_n(x)\big|\leq
    c(r)\delta_n^r(x)\cdot\omega\big(f^{(r)};\delta_n(x)\big).
\end{equation*}

В.~К.~Дзядык (1956) доказал обратную теорему при
$\omega(h)=h^\alpha$ $(\alpha\in(0,1))$. Для этого понадобилось
доказать неравенство для производной $p_n$ типа Маркова--Бернштейна.

Более общую прямую теорему с модулем гладкости второго порядка
$\omega_2$ доказали независимо В.~К.~Дзядык и Г.~Фрейд (G. Freud) в 1959 г.
Третье доказательство см. в \cite[\textbf{5.2.3}]{Timan}.

Одновременное приближение функции и ее производных появилось в
статье автора (1962). А общая прямая теорема с модулем гладкости
$\omega_r$ любого порядка доказана Ю.~А.~Брудным (1974,
\cite{Brudny}, см. \cite{Dzjadyk77}, \cite{DeVore}, \cite{Trigub_Belinsky}).

Аппроксимативную характеристику класса $\gamma W^r[-1,1]$, $r\in\mathbb{N}$, см. в \cite{Trigub73}.

В 60--е годы прошлого века В.К. Дзядык, используя и
развивая методы суммирования рядов Фабера и новые неравенства для
производных многочленов, доказывает прямые и обратные теоремы о
приближении аналитических функций в областях с кусочно гладкой
границей (в его монографии \cite{Dzjadyk77} это стр.~334--490).
Более общие теоремы см. в \cite{Andr_Bel_Dzyad}.

В приведенной выше теореме С.~М.~Никольского (1946) В.~Н.~Темляков (1981) убрал множитель $\ln n$ .

\textbf{Теорема} (\cite{Trigub93}) \emph{Для любого $r\in\mathbb{N}$ существует $c(r)$ такая, что для любой функции $f\in
W^r[-1,1]$ при $n\geq r-1$ найдется последовательность многочленов $p_n$ такая, что при $x\in[-1,1]$}
\begin{equation*}
    \big|f(x)-p_n(x)\big|\leq
    K_r\Big(\frac{\sqrt{1-x^2}}{n+1}\Big)^r+c(r)\frac{(\sqrt{1-x^2})^{r-1}}{(n+1)^{r+1}}.
\end{equation*}

\emph{При этом константу $K_r$ нельзя заменить меньшей, а $c(r)\geq
ce^r$ $(c>0, r\in\mathbb{N})$.}

В.~П.~Моторный \cite{Motorn99} иным методом доказал подобный
асимптотически точный результат при нецелом $r>0$ с константой
Дзядыка и множителем $\ln n$ в остаточном члене, который в его
доказательстве убрать нельзя. Есть еще подобные результаты в
интегральной метрике с весом (см., например, \cite[стр.~248--249]{Trigub_Belinsky}).

См. также подобный результат для приближения того же класса на
полуоси целыми функциями конечной полустепени \cite{Tovstolis}.

Уже изучены многие ограничения в прямых теоремах теории приближений.
Например, односторонние приближения в интегральной метрике уже довольно давно применены к доказательству тауберовых теорем с остаточным членом.

О комонотонных приближениях (функция монотонная и приближающие
полиномы монотонны, функция и полиномы выпуклы и т.д.) написано
много работ (см. гл.~7 в \cite{Dzyadyk_Shevchuk}).

Очевидно, что для $f\in C^r[-1,1]$ $(r\in \mathbb{Z}_+)$
\begin{equation*}
    \begin{gathered}
    f(x)-\sum\limits_{\nu=0}^r\frac{f^{(\nu)}(1)}{\nu!}(x-1)^\nu=\frac{(-1)^{r+1}}{r!}\int\limits_x^1(y-x)^rdf^{(r)}(y)=\int\limits_{-1}^1h_{r,y}(x)df^{(r)}(y),\\
    h_{r,y}(x)=\frac{(x-y)^r}{2(r!)}\big({\rm sign}~(x-y)-1\big).
    \end{gathered}
\end{equation*}

В.К. Дзядык первый хорошо приблизил $h_{r,y}$ (\cite{Dzjadyk77}, гл.~VII, 4], $r=0$). Для комонотонных приближений этот результат использовали
De Vore и Yu (см. \cite{DeVoreYu} или \cite{Dzyadyk_Shevchuk}).

\textbf{Теорема} (\cite{Trigub99}) \emph{Для любых $x\in[-1,1]$ и
$y\in(-1,1)$, а $r$ и $s\in\mathbb{Z}_+$ и $n\geq 2r+1$ существует
многочлен степени не выше $n$ по $x$ такой, что}
\begin{equation*}
\begin{gathered}
    \big|{\rm sign}~(x-y)-p_{n,y}(x)\big|\leq c(r,s)\Big(\frac{1-x^2}{1-x^2+1-y^2+\frac{1}{n^2}|{\rm sign}~x-{\rm
    sign}~y|}\Big)^r\Big(\frac{\delta_n(y)}{|x-y|+\delta_n(y)}\Big)^s,\\
    \delta_n(y)=\frac{\sqrt{1-y^2}}{n}+\frac{1}{n^2}
\end{gathered}
\end{equation*}
(при $r=0$ это неравенство получил В.К.).

Автором изучены прямые теоремы с самыми разными ограничениями:
кусочно односторонние приближения; одновременная аппроксимация с
производными, интерполяцией и учетом положения точки; коэффициенты
многочленов положительные, целые и др. (см. обзор \cite{Trigub2009}).

\bigskip

\textbf{3. Геометрический критерий аналитичности функций.}

\medskip

 Известно, что аналитические (голоморфные) функции на открытом плоском
множестве можно определять по-разному:
$\mathbb{C}$--дифференцируемость, условие Коши--Римана, теорема
Коши--Морера, представимость в окрестности точки степенным рядом,
конформность отображения (два свойства и даже одно).

\textbf{Теорема В.~К.~Дзядыка} \cite{Dzjadyk60} \emph{Рассмотрим в
$\mathbb{R}^3$ три графика (поверхности)}
\begin{equation*}
    z=u(x,y),\quad z=v(x,y),\quad z=\sqrt{u^2(x,y)+v^2(x,y)}.
\end{equation*}

Для того чтобы одна из двух функций $f=u+iv$ и $\overline{f}=u-iv$
была аналитической в области $D$ при непрерывно дифференцируемых $u$
и $v$ необходимо и достаточно, чтобы площади названных трех
поверхностей над любой подобластью $D$ были одинаковыми.

Необходимость проверяется с помощью формулы площади и условий
Коши--Римана. Недавно доказано, что равенства площадей можно
проверять лишь на подобных множествах одного размера (например,
фиксированный многоугольник). Подобное усиление получено и для
теоремы Морера (см. \cite{Volchkov}).

См. также Mathnet.ru

\bigskip

\textbf{4. Некоторые нерешенные вопросы (I--III).}

\medskip

I. Уже давно (1962) получены общие прямые теоремы о приближении
многочленами с целыми коэффициентами на любом отрезке вещественной
прямой длины меньше четырех. См. также \cite{Trigub2001}. Есть и
критерий аппроксимации такими многочленами на множествах
$\mathbb{C}$ (С.~Я.~Альпер, 1964). А вот  теорем типа Дзядыка с целочисленными многочленами нет, если не считать случай квадрата $0\leq {\rm Re}~z, {\rm
Im}~z\leq1$ (\cite{Vit_Volchkov}).

\medskip

II. Автор впервые (1965) построил линейный полиномиальный оператор
$f\mapsto\tau_n$ такой, что при $r\in\mathbb{N}$ и $f\in
C(\mathbb{T})$
\begin{equation*}
    \big\|f-\tau_n(f)\big\|_\infty\asymp\omega_r\Big(f;\frac{1}{n}\Big)\qquad
    (n\in\mathbb{N})
\end{equation*}
(двустороннее неравенство с положительными константами, зависящими лишь от $r$).

При $r=1$ -- это суммы Бернштейна, при $r=2$ -- Рогозинского, а при $r\ge3$ -- специально построенные полиномы типа 
Рогозинского-Бернштейна (см., например, \cite{Dzjadyk77} или \cite{Dzyadyk_Shevchuk}). При $r\geq3$ полиномы 
Бернштейна--Стечкина не подходят для этого, в \cite{Trigub2013} найден специальный модуль гладкости для них.

Указать специальный модуль непрерывности $\omega^*$ такой, что для
любой функции $f\in C(\mathbb{T})$ при некоторой последовательности
$\varepsilon_n\searrow0$, не зависящей от функции,
\begin{equation*}
    \max\limits_x\frac{1}{n+1}\sum\limits_{k=0}^n\big|f(x)-S_n(f;x)\big|\asymp\omega^*(\varepsilon_n).
\end{equation*}

В непериодическом случае подобных результатов с учетом положения
точки (например, с $\delta_n(x)$) нет, если не считать замечательный
результат Тотика (\cite{Totik}, см. также \cite{DeVore}) о
многочленах Бернштейна -- с модулем гладкости Дитциана--Тотика (Z. Ditzian -- V.~Totik).

\medskip

III. В работе \cite{Babenko} исследуется вопрос о поведении множителя
$c(r)$ в теореме Джексона-Стечкина. Тот же вопрос о теоремах А.~Ф.~Тимана и
Ю.~А.~Брудного (см. в п.~2). И в оценке приближения снизу, конечно. 

\end{document}